\documentclass[11pt]{article}
\usepackage{amsmath,amsfonts,latexsym,amssymb,amscd, graphicx,pictexwd, mathrsfs, dsfont, color,mathtools}
\newcommand{\rb}{{\mathrm b}}
\newcommand{\rh}{{\mathrm h}}

\newcommand{\rnw}{{\mathrm{nw}}}

\renewcommand{\P}{{\mathbb P}}
\newcommand{\E}{{\mathbb E}\,}

\newcommand{\R}{{\mathbb R}}

\newcommand{\ZZ}{{\mathbb Z}}
\newcommand{\pp}{{\mathbf p}}

\newcommand{\cc}{{\mathbf c}}
\newcommand{\zz}{{\mathbf z}}

\newcommand{\xx}{{\mathbf x}}
\newcommand{\uu}{{\mathbf u}}
\newcommand{\vv}{{\mathbf v}}
\newcommand{\yy}{{\mathbf y}}
\newcommand{\dd}{{\mathbf d}}
\newcommand{\0}{{\mathbf 0}}

\newcommand{\cA}{{\mathcal A}}

\newcommand{\fh}{{\mathfrak h}}
\newcommand{\fb}{{\mathfrak b}}
\newcommand{\fg}{{\mathfrak g}}

\newcommand{\fd}{{\mathfrak d}}

\newcommand{\sfC}{{\mathsf{C}}}
\newcommand{\sfUC}{{\mathsf{UC}}}

\newcommand{\Exp}{{\rm Exp}}

\newtheorem{thm}{Theorem}

\newtheorem{coro}{Corollary}[thm]
\newtheorem{lem}{Lemma}[section]
\newtheorem{prop}{Proposition}
\newtheorem{rem}{Remark}

\numberwithin{equation}{section}

\begin{document}
\title{Ergodicity of the KPZ Fixed Point}
\author{Leandro P. R. Pimentel}
\date{\today}
\maketitle

\begin{abstract}  
In this paper we consider the Kardar-Parisi-Zhang (KPZ) fixed point $\left(\fh_t\,,\,t\geq 0\right)$ \cite{CoQuRe,MaQuRe} and prove that, for suitable initial conditions, $\fh_t(x)-\fh_t(0)$ converges to a two-sided Brownian motion with zero drift and diffusion coefficient $2$, as $t\to\infty$. The heart of the proof is the coupling method, that allows us to compare the TASEP height function started from a perturbation of density $1/2$ with its invariant counterpart which, under KPZ scaling, turns into uniform estimates for the KPZ fixed point. 
\end{abstract}

\section{Introduction and Main Result}
 
In $d+1$ stochastic growth models the object of interest is a height function, $h(x,t)$ for $x\in\R^d$ (space) and $t\geq 0$ (time),  whose evolution is described by a random mechanism. For fairly general growth models one has a deterministic macroscopic shape for the $h(x,t)$ and its fluctuations, under proper space-time scaling, are expected to be characterized by a universal distribution. The Kardar-Parisi-Zhang (KPZ) fixed point is a Markov process introduced by Corwin, Quastel and Remenik \cite{CoQuRe} and Matetski, Quastel and Remenik  \cite{MaQuRe} that gives the limit fluctuations of a wide class of growth models with $d=1$ that possesses a local slope dependent growth rate and a smoothing mechanism, combined with space-time random forcing with rapid decay of correlations.  The KPZ equation \cite{KPZ}
$$\frac{\partial h}{\partial t}=\frac{1}{2}\left(\frac{\partial h}{\partial x}\right)^2+\frac{\partial^2 h}{\partial^2 x }+\xi\,,$$
where $\xi$ is a space-time white noise, is a canonical example of such a growth model, providing its name to the universality class \cite{AmCoQu,BoCoFeVe}. It is conjectured that for such models
$$h(an^{2/3}x,nt)\stackrel{dist.}{\sim}btn+cn^{1/3}\fh_t(x)\,,$$
for some model dependent constants $a,b,c\in\R$, where $\fh_t(x)$ is a universal space-time process, called the KPZ-fixed point.
\newline

Another canonical example of KPZ fluctuations is given by the height function associated to the totally asymmetric simple exclusion process (TASEP), where the limiting behaviour was stablished initially by Johansson \cite{Jo1,Jo2}. In the past twenty years there has been a significant amount of improvements of the theory. The exact statistics for different initial interface profiles, resulting in different types of Airy processes, were computed using the TASEP height function and its connection with integrable probability. Recently, a unifying approach based on a previous work of Sasamoto and coauthors \cite{BoFePrSa,Sa} was developed by Matetski, Quastel and Remenik \cite{MaQuRe}. They derived a Fredholm determinant formula for the finite dimensional distributions of the TASEP height function, and proved that this formula conducted to the computation of the transition probabilities of the KPZ fixed point.
\newline

The evolution of the KPZ fixed point takes place on the space $\sfUC$ of upper semicontinuous (generalized) functions $\fh:\R\to [-\infty,\infty)$ such that $\fh(x)<C(1+|x|)$ for some $C>0$ (we allow the initial data to take value $-\infty$). Recall that $\fh$ is upper semicontinuous if and only if its hypograph $\{(x,y)\,:\,y\leq \fh(x)\}$ is closed in $[-\infty,\infty)\times \R$. The topology of local convergence in $\sfUC$ can be defined in terms of the Hausdorff distance between hypographs (see Section 3.1 in \cite{MaQuRe}).  Let $\fh_t\equiv\fh_t(\cdot;\fh)$ denote the KPZ fixed point starting at $\fh_0=\fh$. For a fixed time $t>0$ the distribution of $\fh_t(\cdot;\fh)$ is determined by       
\begin{equation}\label{Int}
\P\Big(\fh_t(x;\fh)\leq \fg(x)\,,\,\forall\,x\in\R\Big)=\det\left(I-K_{t,\fh,\fg}\right)_{L^2(\R)}\,,
\end{equation}
where $\fg$ is lower semicontinuous function, i.e. $-\fg\in\sfUC$. For a description of the operator $K_{t,\fh,\fg}$ in terms of $\fh$ and $\fg$ the author address the reader to \cite{MaQuRe}. The determinant on the right hand side is a Fredholm determinant on the Hilbert space $L^2(\R)$, and $I$ is the identity:
$$\det(I-K)_{L^2(\R)}:=\sum_{n=0}^{\infty}\frac{(-1)^n}{n!}\int_{\R^n}\det\left(K(x_i,x_j)\right)_{i,j=1}^n dx_1\cdots dx_n\,,$$ 
where $K(x,y)$ is the kernel of $K$. From this formula one can recover several of the classical Airy processes by starting with special profiles for which the respective operators $K$ are explicit (see Section 4.4 of \cite{MaQuRe}). For instance, the Airy$_2$ process $\cA(\cdot)=\fh(\cdot;\fd)$ is defined by taking the initial profile $\fh=\fd$ where $\fd(0)=0$ and $\fd(z)=-\infty$ for all $z\neq 0$.
\newline

The KPZ fixed point satisfies some fundamental symmetries are expressed below for fixed $t>0$ in terms of distributional equalities as processes in $\sfUC$ \cite{MaQuRe}.
\begin{itemize}
\item[(i)] 1-2-3 Scaling: $S_{\gamma^{-1}}\fh_{\gamma^{-3}t}(\cdot;S_\gamma\fh)\stackrel{dist.}{=}\fh_{t}(\cdot;\fh)$, where $S_\gamma\fh(x):=\gamma^{-1}\fh(\gamma^{2} x)$ for $\gamma>0$.
\item[(ii)] Affine Symmetry: $\fh_t(\cdot;\Upsilon_{\alpha,\beta}\fh)\stackrel{dist.}{=}\Upsilon_{\alpha,\beta+\frac{\alpha^2}{4}t}T_{\frac{\alpha}{2} t}\fh_t\left(\cdot;\fh\right)$, where 
$$\Upsilon_{\alpha,\beta}\fh(x):=\beta+\alpha x +\fh(x)\,\,\mbox{ and  }\,\,T_z \fh(x):=\fh(z+x) \,.$$ 
\item[(iii)] Space Stationarity: $T_z\fh_t\left(\cdot;T_{-z}\fh\right)\stackrel{dist.}{=}\fh_t(\cdot,\fh)$. 
\end{itemize}
We note that, by the 1-2-3 scaling (i) with $\gamma=t^{1/3}>0$, we have that  
\begin{equation}\label{123}
\fh_t(\cdot;\fh)\stackrel{dist.}{=}S_{t^{-1/3}}\fh_1(\cdot;S_{t^{1/3}}\fh)\,.
\end{equation}
Another important symmetry is related to time invariance: let $\fb$ be a two-sided Brownian motion with diffusion coefficient $2$, then 
$$\Delta \fh_t(\cdot;\fb)\stackrel{dist.}{=}\fb(\cdot)\,\mbox{ for all }\,t\geq 0\,,\,\mbox{ where }\Delta\fh(x):=\fh(x)-\fh(0)\,.$$ 
Combining this with (ii) (affine invariance), one sees that the drifted Brownian motion $\Upsilon_{0,\alpha}\fb$ is also invariant. The aim of this article is to prove convergence of $\Delta \fh_t$ to $\fb$, as $t\to\infty$, for suitable initial conditions (ergodicity). 
\newline

The study of local space regularity in the KPZ universality class was initiated by H\"agg \cite{Hag}, where it was proved local finite-dimensional convergence of the Airy$_2$ process to Brownian motion, using asymptotic estimates related to a determinantal structure similar to \eqref{Int}. This result was extended to functional convergence in \cite{CaPi}, by means of the coupling method. Using kernel estimates for the operator in \eqref{Int}, Matetski, Quastel and Remenik \cite{MaQuRe} proved that for every $\fh\in \sfUC$ and $t>0$, $\fh_t(\cdot;\fh)$ has H\"older $1/2-$ regularity in space, and that 
$$S_{\sqrt{\epsilon}}\Delta\fh_t(x)=\epsilon^{-1/2}\Big(h_t(\epsilon x;\fh)-h_t(0;\fh)\Big)\,$$
converges to $\fb$, as $\epsilon\to 0$, in terms of finite dimensional distributions. Functional convergence was proved by Pimentel \cite{Pi} for several versions of Airy processes, and stronger forms of local Brownian behaviour were proved by Corwin and Hammond \cite{CoHa} and Hammond \cite{Ha}. From the 1-2-3 scaling \eqref{123}, the long time behaviour of $\Delta\fh_t$ can be rewriten in terms of the local behaviour of the process at time $t=1$,
$$\Delta\fh_t(\cdot;\fh) \stackrel{dist.}{=} S_{\sqrt{\epsilon}}\Delta\fh_1(\cdot;S_{\sqrt{\epsilon^{-1}}}\fh)\,,\,\mbox{ with }\,\epsilon=t^{-2/3}\,,$$
which allows one to obtain ergodicity (in terms of finite dimensional distributions) from the local convergence to Brownian motion with zero drift, as soon as $S_{\sqrt{\epsilon^{-1}}}\fh$ has a limit. 
\newline

Recently there has been considerable developments in describing the space-time structure of the KPZ fixed in terms of a variational formula \cite{CoLiWa,DaOrVi,FeOc,MaQuRe} involving the Airy$_2$ process: for fixed $x\in\R$,
\begin{equation}\label{VarFor}
\fh_t(x;\fh)\stackrel{dist.}{=}\sup_{z\in\R}\left\{\fh(z)+t^{1/3}\cA(zt^{-2/3})-\frac{(z-x)^2}{t}\right\}\,.
\end{equation}     
(See Theorem 4.18 and Remark 4.19 in \cite{MaQuRe}.) In this paper we rely on \eqref{VarFor} and use the coupling method applied to a last-passage percolation discrete approximation of $\fh_t$ to prove long time convergence in probability. 

\begin{thm}\label{Coupling}
Assume there exist $\gamma_0>0$ and $\psi(r)$ such that for all $\gamma>\gamma_0$ and $r\geq 1$
\begin{equation}\label{Assump}
\P\left(S_\gamma\fh(x)\leq r|x|\,,\,\forall\,|x|\geq 1\right)\geq 1-\psi(r)\,,\mbox{ where }\lim_{r\to\infty}\psi(r)=0\,.
\end{equation}
Let $a>0$ and $t\geq a^{3/2}$. There exists a coupling $\Big(\fh_t(\cdot;\fh),\fh_t(\cdot;\fb)\Big)$, where $\fh$ and $\fb$ are sample independently at $t=0$, and a function $\theta(\delta)$ such that  
\begin{equation}\label{Coupling1}
\P\left(\sup_{x\in[-a,a]}|\Delta\fh_t(x;\fh)-\Delta\fh_t(x;\fb)|>\eta \sqrt{a}\right)\leq \theta(at^{-2/3})+\frac{(at^{-2/3})^{1/4}}{\eta}\,,\,\forall\,\eta>0\,,
\end{equation}
where $\lim_{\delta\to 0}\theta(\delta)=0$.
\end{thm}

\begin{rem}\label{Transversal}
The function $\theta$ might depend on $\psi$, but not on the parameters $a,t>0$. Its behaviour is also  connected to transversal fluctuations of a maximal path in the (exponential) last-passage percolation model under the $2/3$ scaling (Proposition \ref{InterControl2}).  
\end{rem}

As a consequence of Theorem \ref{Coupling}, one can let $a=a_t$ (e.g. $a_t=t^{\kappa}$ for $\kappa\in(0,2/3)$) and still get convergence to equilibrium under proper rescaling.
\begin{coro}
Under \eqref{Assump}, if $\lim_{t\to\infty}a_{t} t^{-2/3}=0$ then 
$$\lim_{t\to\infty}\P\left(\sup_{x\in[-a_t,a_t]}|\Delta\fh_t(x;\fh)-\Delta\fh_t(x;\fb)|>\eta \sqrt{a_t}\right)=0\,,\,\,\mbox{ as }t\to\infty\,.$$
In particular, since $\Delta\fb_t$ is invariant under diffusive scaling, 
$$\lim_{t\to\infty}S_{\sqrt{a_t}}\Delta\fh_t\stackrel{dist.}{=}\fb\,.$$
\end{coro}

\noindent\paragraph{\bf Acknowledgements} The author would like to thank Patrik Ferrari and Daniel Remenik for useful comments and  enlightening discussions concerning a previous version of this article. This research was supported in part by the International Centre for Theoretical Sciences (ICTS) during a visit for participating in the program - Universality in random structures: Interfaces, Matrices, Sandpiles (Code: ICTS/urs2019/01), and also by the CNPQ grants 421383/2016-0 and 302830/2016-2, and by the FAPERJ grant E-26/203.048/2016.   

\section{The TASEP Growth Process and Last-Passage Percolation}\label{TASEPLPP}

The totally asymmetric simple exclusion process (TASEP) is a Markov process $\left(\eta_t\,,\,t\geq 0\,\right)$ with state space $\{0,1\}^\ZZ$. When $\eta_t(k)=1$, we say that  site $x$ is occupied by a particle at time $t$, and it is empty (hole) if $\eta_t(k)=0$.  Particles jump to the neighbouring right site with rate $1$ provided that the site is empty (the exclusion rule). For a initial $\eta\in\{0,1\}^\ZZ$, let $N^\eta_t$ denote the total number of particles which jumped from site $0$ to site $1$ during the time interval $[0,t]$. The TASEP growth model is represented by a height function  $h_t(\cdot)=h_t(\cdot;\eta)\in\ZZ^\ZZ$, defined for each $t\geq 0$ as 
$$h_t(k;\eta)=\left\{\begin{array}{lll}
2N^\eta_t+\sum_{j=1}^{k}(1-2\eta_t(j)) & \mbox{ for } k\geq 1\\
2N^\eta_t & \mbox{ for } k=0\\
2N^\eta_t-\sum_{j=k+1}^{0}(1-2\eta_t(j)) &\mbox{ for } k\leq -1\,.\end{array}\right.$$ 
We extend $h_t$ to a function $h_t:\R\to\R$ by linear interpolation of its values at integer numbers. The time evolution of $h_t$, induced by the TASEP dynamics, is that local minima become local maxima at rate $1$: if $h_t(z\pm 1)=h_t(z)+1$ then $h_t(z)\to h_t(z)+2$ at rate $1$. For a given sequence of initial profiles $\left(\eta^{(n)}\,,\,n\geq 1\right)$, let\footnote{In \cite{MaQuRe} they considered the growth interface moving downwards and the scaling factor $\epsilon\to 0$, while we are considering it moving upwards and setting $\epsilon=n^{-2/3}$.}
\begin{equation*}\label{HeightScaling}
\fh_{t,n}\left(x;\eta^{(n)}\right):=\frac{tn-h_{2tn}\left(2xn^{2/3};\eta^{(n)}\right)}{n^{1/3}}\,,\,\mbox{ for $x\in\R$}\,.
\end{equation*}
Matetski, Quastel and Remenik \cite{MaQuRe} proved that if
\begin{equation}\label{kpz1c}
\lim_{n\to\infty}\fh_{0,n}(\cdot;\eta^{(n)})\stackrel{\sfUC}{=}\fh(\cdot)\,,
\end{equation}
then
\begin{equation}\label{kpz1}
\lim_{n\to\infty}\fh_{t,n}(\cdot;\eta^{(n)})\stackrel{dist.}{=}\fh_t(\cdot;\fh)\,.
\end{equation}
\newline
Assumption \eqref{kpz1c} restricts the study of the scaling limit to perturbations of density $1/2$, where one has a limit at time zero under diffusive scaling. We note that for any $\fh\in \sfUC$, with $\fh(0)=0$, one can construct a sequence of initial particle configurations $\eta^{(n)}$ such that $\fh_{0,n}(\cdot)=\fh_{0,n}(\cdot;\eta^{(n)})\stackrel{\sfUC}{\to}\fh(\cdot)$.
\newline

The TASEP growth process has a equivalent formulation in terms of a corner growth model. Let $\rh(\eta)=(\rh(k;\eta))_{k\in\ZZ}\in(\ZZ^2)^\ZZ$ denote the down-right nearest-neighbour path in $\ZZ^2$ constructed from the particle configuration $\eta$ as follows: 
$$\rh(0;\eta)=(0,0)\,\mbox{ and }\,\rh(k+1;\eta)=\rh(k;\eta)+\eta(k+1)(0,-1)+(1-\eta(k+1))(1,0)\,.$$ 
The path $\rh=\rh(\eta)$ splits $\ZZ^2$ into two regions and we denote $\Gamma$ the one that includes $\rh$ and the negative axes. We will also assume that $\eta$ has a positive density of particles to the left of the origin, and a positive density of holes to the right of the origin. The corner growth model $(\Gamma_t)_{t\geq 0}$ is described by the set $\Gamma_t$ of occupied vertices at time $t$, where $\Gamma_0=\Gamma$ and each site $(k,l)\in\Gamma_t^c$ becomes occupied at rate $1$, once the sites $(k-1,l)$ and $(k,l-1)$ are both occupied. For each $(k,l)\in\Gamma$ set $L^{\rh}(k,l):=0$, and for $(k,l)\in\Gamma^c$ let 
$$L^{\rh}(k,l):=\mbox{ the time when $(k,l)$ becomes occupied}\,,$$
so that 
$$\Gamma_t=\left\{(k,l)\in\ZZ^2\,:\,L^{\rh}(k,l)\leq t\right\}\,.$$
It is enough clear that the boundary of $\Gamma_t$ and the TASEP height function $h_t$ have similar evolution rules (up to a $45^o$ rotation). Thus, if we construct $\rh$ using the same initial particle configuration $\eta$, and match the transition rates, then 
\begin{equation}\label{TALPP}
L^{\rh}(k,l)\leq t \,\,\Longleftrightarrow\,\, h_t(k-l)\geq k+l\,.
\end{equation}

By combining \eqref{TALPP} together with \eqref{kpz1}, the KPZ fixed point can be obtained as the limit fluctuations of $L^\rh$ as well. Denote 
$$[k]_{n}\equiv (n+k,n-k)\,\mbox{ and }\,L^{\rh}[k]_n\equiv L^{\rh}\left([k]_n\right)\,.$$  
For $x\in\R$ with $|x|\leq n$ we set $L^{\rh}[x]_{n}$ by linear interpolation of the values of $L^{\rh}[\cdot]_{n}$ at integer numbers. For a given sequence of initial profiles $\left(\eta^{(n)}\,,\,n\geq 1\right)$ denote $\rh^{(n)}=\rh\left(\eta^{(n)}\right)$ and define the process
$$H_{t,n}\left(x\right)\equiv H_{t,n}\left(x;\rh^{(n)}\right) = \frac{L^{\rh^{(n)}}[2^{2/3} xn^{2/3}]_{\lfloor tn\rfloor}-4tn}{2^{4/3} n^{1/3}}\,,$$
 if $|2^{2/3} xn^{2/3}|\leq n$, and set $H_{t,n}(x)= 0$ otherwise. In view of \eqref{TALPP}, convergence of the finite dimensional distributions of $H_{t,n}$ to the finite dimensional distributions of $\fh_t$ follows from \eqref{kpz1}. It is also known that $\left\{H_{t,n}\,:\,n\geq 1\right\}$ is relatively compact \cite{CoLiWa,FeOc,Pi} with respect to weak convergence of measures in the space of continuous real functions on compact sets, endowed with the uniform norm \cite{Bi}, which implies functional convergence of $H_{t,n}$ to the KPZ fixed point: under \eqref{kpz1c}, 
\begin{equation}\label{kpz3}
\lim_{n\to\infty}H_{t,n}(\cdot;\rh^{(n)})\stackrel{dist.}{=}\fh_t(\cdot;\fh)\mbox{ for all }t\geq 0\,.
\end{equation}

\subsection{Last-Passage Percolation and Maximal Paths}\label{MaxPath}
Recall that $L^{\rh}(k,l)$ denote the time when $(k,l)$ becomes occupied, so that the corner growth model can be expressed as 
$\Gamma_t=\left\{(k,l)\in\ZZ^2\,:\,L^{\rh}(k,l)\leq t\right\}$.
Since each site $(k,l)\in\Gamma_t^c$ becomes occupied at rate $1$, once the sites $(k-1,l)$ and $(k,l-1)$ are both occupied, we have that 
\begin{equation}\label{rec}
L^{\rh}(k,l)=\omega_{k,l}+\max\{L^{\rh}(k-1,l)\,,\,L^{\rh}(k,l-1)\}\,,
\end{equation}
where the random variables $\omega_{k,l}$ are i.i.d. with an exponential distribution of parameter $1$. This allows us to define the corner growth process in terms of a last-passage percolation model (LPP) as follows. For $\xx,\yy\in\ZZ^2$ with $\xx \leq \yy$ (coordinate-wise), we say that a sequence $\pi=(\xx_0,\xx_1,\cdots,\xx_j)$ is an up-right path from $\xx$ to $\yy$, if $\xx_{i+1}-\xx_{i}\in\{(1,0),(0,1)\}$, $\xx_0=\xx$ and $\xx_j=\yy$. We denote $\Pi(\xx,\yy)$ the set compose of all up-right paths from $\xx$ to $\yy$. The random environment in our setting is given by a collection $\omega\equiv\left\{\omega_{k,l}\,:\,\in\ZZ^2\right\}$ of i.i.d. random variables (passage times) with exponential distribution of parameter $1$. The (point-to-point) last-passage percolation time is defined as
$$L(\xx,\yy):=\max_{\pi\in\Pi(\xx,\yy)}\sum_{\zz\in\pi} \omega_{\zz}\,.$$
We do not have super-additivity for $L(\xx,\yy)$, since for each path $\pi\in\Pi(\xx,\yy)$ we include both end points in the sum of the passage times. This implies that, when concatenating paths, we sum twice the passage time at the concatenation point.  However, it holds that
\begin{equation}\label{super}
L(\xx,\yy)\geq L(\xx,\zz)+L(\zz,\yy)-\omega_\zz\,,
\end{equation}    
for all $\xx\leq \zz\leq \yy$.
\newline

To establish the connection with the corner growth model we need to define curve-to-point last-passage percolation times. In order to do so, we say that $\zz\in\ZZ^2$ is a concave corner of $\rh$ if $\zz-(1,1)=\rh(k)$ for some $k\in\ZZ$ and  
$$\rh(k-1)=\rh(k)-(0,1)\,\mbox{ and }\,\rh(k+1)=\rh(k)+(1,0)\,.$$
For each $\xx\in\Gamma^c$, we denote $\sfC_\xx^{\rh}$ the set of all concave corners $\zz$ of $\rh$ such that $\zz\leq \xx$. Thus, \eqref{rec} implies that for all $\xx\in\Gamma^c$,
$$L^{\rh}(\xx):=\max_{\zz\in\sfC_\xx^{\rh}} L(\zz,\xx)\,.$$ 
\newline

The maximal path from $\rh$ to $\xx\in\Gamma^c$ is the a.s. unique up-right path $\pi^{\rh}(\xx)$, starting at some corner in $\sfC_\xx^{\rh}$ end ending at $\xx$ such that 
$$L^{\rh}(\xx)=\sum_{\zz\in\pi^{\rh}(\xx)}\omega_\zz\,.$$ 
Given a maximal path $\pi^{\rh}(\xx)$, with $\xx=(k,l)> 0$ (coordinate-wise), consider  the last point of $\pi^{\rh}(\xx)$ (in the up-right orientation) that intersects the positive coordinate axis $\{\zz=(z,0)\mbox{ or }\zz=(0,z)\mbox{ and }z>0\}$. If such a point does not exist we must have that $\pi^{\rh}(\xx)$ starts at $(1,1)$ (recall that $\rh$ is a down-right path and that $\rh(0)=(0,0)$), and in such a case we define $Z=1$. To distinguish between intersections via the horizontal- and vertical-axis we introduce a non-zero integer-valued random variable $Z=Z^{\rh}(\yy)$ such that if $Z>0$ then the intersection point is $(Z,0)$, while if $Z<0$ then the intersection point is $(0,-Z)$. \newline

\subsection{Last-Passage Percolation with Boundary Condition}  
To apply the coupling method, it will be more convenient to work with a LPP model with boundary condition instead of having a down-right path as an initial profile. The boundary condition 
$$\rb:=\{\omega^\rb_\xx\,:\,\mbox{ for }\xx=(z,0)\mbox{ or }\xx=(0,z)\mbox{ and }z\geq 0\}\,,$$
is given by real numbers $\omega_{\xx}^\rb\geq 0$ that are placed along the non-negative coordinate axis, and we will always assume that $\omega^\rb_{(0,0)}=0$. Define
$$\rb(z)=\left\{\begin{array}{lll}
\sum_{i=1}^{-z}\omega^\rb_{0,i} & \mbox{ for } z<0\\
0 & \mbox{ for } z=0\\
\sum_{i=1}^z\omega^\rb_{i,0} &\mbox{ for } z>0\,,\end{array}\right.$$
and for $(k,l)>(0,0)$ and $z\in[-l,k]$, 
$$L_z(k,l)=\left\{\begin{array}{lll}
L\left((1,-z),(k,l)\right) & \mbox{ for } z\in[-l,0)\\
L\left((z,1),(k,l)\right) &\mbox{ for } z\in(0,k]\,.\end{array}\right.$$
The last-passage percolation time to $\xx=(k,l)>(0,0)$, with boundary condition $\rb$, is defined as 
\begin{equation}\label{LPP}
\bar L^{\rb}(k,l):=\max_{z\in [-l,k]\setminus\{0\}} \left\{ \rb(z) + L_z(k,l)\right\}\,.
\end{equation}
The exit point associated to the boundary condition $\rb$ is defined as 
$$\bar Z^{\rb}(k,l) :=\max\arg\max_{z\in [-l,k]\setminus\{0\}} \left\{ \rb(z) + L_z(k,l)\right\}\,,$$
which corresponds to the right-most point $z\in[-l,k]\setminus\{0\}$ for which $\rb(z) + L_z(k,l)=\bar L^{\rb}(k,l)$. It is clear that for fixed $n\geq 1$, 
\begin{equation}\label{IncExit}
\bar Z^{\rb}\left(i+n,i-n\right) \,\mbox{ is a non decreasing function of $i\in[-n,n]$\,.}
\end{equation}

\noindent\paragraph{\bf Example 1: Invariant Regime.} 
Let $\rho\in(0,1)$ and define a collection of independent random variables with $\omega^\rho_{(0,0)}=0$,   
$$\omega_{(z,0)}^\rho\stackrel{dist.}{=}\Exp_1(1-\rho)\,\mbox{ and }\,\omega_{(0,z)}^\rho\stackrel{dist.}{=}\Exp_1(\rho)\,\mbox{ for }z\geq 1\,.$$
This boundary condition has a connection with the TASEP process started from Bernoulli($\rho$) measure conditioned on $\eta(0)=0$ (hole at the origin) and $\eta(1)=1$ (particle at site one) \cite{BaCaSe}. A fundamental property that we will use is that the increments of the last-passage times $\bar L^{\rho}$ along the anti-diagonal are i.i.d. with a well known distribution.  Precisely, define (recall that $[k]_n:=(n+k,n-k)$)
\begin{equation}\label{statinc}
\zeta^{\rho}_{k,n}:=\bar L^\rho[k]_n-\bar L^\rho[k-1]_n\,,\mbox{ for }\,k=-n+1,\dots,n\,.
\end{equation}
Then $\left\{\zeta^\rho_{k,n}\,:\,k=-n+1,\dots,n\right\}$ is a collection of independent random variables with 
\begin{equation}\label{stat}
\zeta^\rho_{k,n}\stackrel{dist.}{=}\Exp_1(1-\rho)-\Exp_2(\rho)\,,
\end{equation}
where $\Exp_1(1-\rho)$ and $\Exp_2(\rho)$ are independent random variables with exponential distribution of parameter $\rho$ and $1-\rho$, respectively (see Lemma 4.2 in \cite{BaCaSe}). 

\noindent\paragraph{\bf Example 2: Boundary induced by Curve to Point LPP.}
Let $\omega^{\rh}(0,0)=0$, and for $z\geq 1$ set 
$$\omega^{\rh}(z,0):=L^{\rh}(0,z)-L^{\rh}(z-1,0)\,\mbox{ and }\,\omega^{\rh}(0,z):=L^{\rh}(0,z)-L^{\rh}(0,z-1)\,.$$
(Notice that $L^{\rh}(0,0)=0$ since $(0,0)\in\Gamma$.) Then 
\begin{equation}\label{Recover}
\bar L^{\rh}(k,l)=L^{\rh}(k,l)\,\mbox{ and }\,\bar Z^{\rh}(k,l)=Z^{\rh}(k,l)\,,\mbox{ for all }k,l>0\,.
\end{equation}
To see this, first consider $z\in(0,k]$. Then,  
$$\rb^{\rh}(z)+L_{z}(k,l)=L^{\rh}(0,z)+L_{z}(k,l)\leq L^{\rh}(k,l)\,.$$  
Similarly we have the same inequality for $z\in[-l,0)$, and hence $\bar L^\rh(k,l)\leq L^{\rh}(k,l)$. On the other hand, if one take $Z=Z^\rh(k,l)$ (as defined at the end of Section \ref{MaxPath}) then  
$$L^{\rh}(k,l)=\rb^{\rh}(Z)+L_{Z}(k,l)\leq \bar L^{\rh}(k,l)\,.$$ 
and thus $L^\rh(k,l)\leq \bar L^{\rh}(k,l)$.
\newline
  
\subsection{Argmax Comparison and Attractiveness under the Basic Coupling}\label{BasCoup}
  
Given $\rb$ and $\tilde\rb$, the basic coupling is the joint construction of $(L^{\rb},L^{\tilde\rb})$ using the same point-to-point last-passage percolation times \eqref{LPP}:
\begin{equation}\label{basiccoupling}
\bar L^{\rb}(k,l):=\max_{z\in [-l,k]\setminus\{0\}} \left\{ \rb(z) + L_z(k,l)\right\}\,\mbox{ and }\,\bar L^{\tilde\rb}(k,l):=\max_{z\in [-l,k]\setminus\{0\}} \left\{ \tilde\rb(z) + L_z(k,l)\right\}\,.
\end{equation}
The next lemma allow us to compare, under the basic coupling, local increments associated to different boundary conditions by looking at the relative positions of the respective exit points.  
    
\begin{lem}[Argmax Comparison]\label{comparison}
Under \eqref{basiccoupling},
$$\mbox{  if }\,\,\bar Z^{\rb}[j]_n\leq \bar Z^{\tilde\rb}[i]_n\,,\,\mbox{ for }\,i\leq j\,,\,\mbox{  then  }\,\,\bar L^{\rb}[j]_n-\bar L^{\rb}[i]_n \leq \bar L^{\tilde\rb}[j]_n-\bar L^{\tilde\rb}[i]_n\,.$$
\end{lem}

\noindent{\bf Proof\,\,} This lemma is similar to Lemma 2.1 in \cite{Pi}, and the proof follows the same lines. Denote $\pi_z(\xx)$ the maximal path associated to $L_{z}(\xx)$. Denote $z_1\equiv \bar Z^{\rb}[j]_n$ and $z_2\equiv \bar Z^{\tilde\rb}[i]_n$. Let $\cc$ be a crossing between $\pi_{z_1}([j]_n)$ and $\pi_{z_2}([i]_n)$. Such a crossing always exists because, by assumption, $i\leq j$ and $z_1\leq z_2$. We remark that, by \eqref{super},
$$\bar L^{\tilde\rb}[j]_n  \geq  \tilde\rb(z_2) + L_{z_2}([j]_n) \geq  \tilde\rb(z_2) + L_{z_2}(\cc) + L(\cc,[j]_n)-\omega_\cc \,.$$
We use this, and that (since $\cc\in\pi_{z_2}([i]_n)$)
$$ \tilde\rb(z_2) + L_{z_2}(\cc)-\bar L^{\tilde\rb} [i]_n= -L(\cc,[i]_n)+\omega_\cc\,,$$
in the following inequality:
\begin{eqnarray*}
\bar L^{\tilde\rb}[j]_n - \bar L^{\tilde\rb}[i]_n & \geq & \tilde\rb(z_2)+L_{z_2}(\cc) + L(\cc,[j]_n)-\omega_\cc - L^{\tilde\rb}[i]_n\\
& = & L(\cc,[j]_n) - L(\cc,[i]_n)\,.
\end{eqnarray*}
By \eqref{super},
$$ - L(\cc\,,\,[i]_n)\geq \bar L^{\rb}(\cc)-\bar L^{\rb}[i]_n-\omega_\cc\,,$$
and hence (since $\cc\in\pi_{z_1}([j]_n)$)
\begin{eqnarray*}
\bar L^{\tilde\rb}[j]_n - \bar L^{\tilde\rb}[i]_n & \geq & L(\cc,[j]_n) - L(\cc,[i]_n)\\
& \geq & L(\cc,[j]_n) + \bar L^{\rb}(\cc)-L^{\rb}[i]_n-\omega_\cc\\
& = & \bar L^{\rb}(\cc)+\left(L(\cc,[j]_n) -\omega_\cc\right)-L^{\rb}[i]_n\\
& = & \bar L^{\rb}[j]_n-\bar L^{\rb}[i]_n\,.
\end{eqnarray*}

\hfill$\Box$\\ 

Another useful property is attractiveness of the LPP model with boundary condition. It states that, under the basic coupling, if one starts the with ordered boundary conditions then the last-passage percolation times remain ordered as well.     

\begin{lem}[Attractiveness]\label{attract}
Under  \eqref{basiccoupling}, if $\rb(j)-\rb(i)\leq \tilde\rb(j)-\tilde\rb(i)$, for all $i\leq j$, then  
$$\bar L^{\rb}[j]_n-\bar L^{\rb}[i]_n\leq \bar L^{\tilde\rb}[j]_n-\bar L^{\tilde\rb}[i]_n\,,\,\forall\,i\leq j\,.$$
\end{lem}

\noindent{\bf Proof\,\,} Denote again $z_1\equiv Z^{\rb}[j]_n$ and $z_2:=Z^{\tilde\rb}[i]_n$. If $z_1 \leq z_2$ then it follows from Lemma \ref{comparison} (we do not need to use the assumption). If  $z_1>z_2$ then
\begin{eqnarray*}
\bar L^{\tilde\rb}[j]_n-\bar L^{\tilde\rb}[i]_n &-&\left(\bar L^{\rb}[j]_n-\bar L^{\rb}[i]_n\right)\\
&=&\bar L^{\tilde\rb}[j]_n-\big(\tilde\rb(z_2)+L_{z_2}[i]_n\big)-\Big(\big(\rb(z_1)+L_{z_1}[j]_n\big)-\bar L^{\rb}[i]_n\Big)\\
&=&\bar L^{\tilde\rb}[j]_n-\big(\tilde\rb(z_2)+L_{z_1}[j]_n\big)-\Big(\big(\rb(z_1)+L_{z_2}[i]_n\big)-\bar L^{\rb}[i]_n\Big)\\
&=&\bar L^{\tilde\rb}[j]_n-\big(\tilde\rb(z_2)+L_{z_1}[j]_n\big)+\Big (\bar L^{\rb}[i]_n- \big(\rb(z_1)+L_{z_2}[i]_n\big)\Big)\\
&=&\bar L^{\tilde\rb}[j]_n-\big(\tilde\rb(z_1)+L_{z_1}[j]_n\big)+\Big(\bar L^{\rb}[i]_n- \big(\rb(z_2)+L_{z_2}[i]_n\big)\Big)\\
&+&\big(\tilde\rb(z_1)-\tilde\rb(z_2)\big)-\big(\rb(z_1)-\rb(z_2)\big) \,.
\end{eqnarray*}
On the other and,
$$\bar L^{\tilde\rb}[j]_n-\big(\tilde\rb(z_1)+L_{z_1}[j]_n\big)\geq 0\,, $$
and 
$$\bar L^{\rb}[i]_n- \big(\rb(z_2)+L_{z_2}[i]_n\big)\geq 0\,,$$
while, by assumption, 
$$\tilde\rb(z_1)-\tilde\rb(z_2)\geq \rb(z_1)-\rb(z_2)\,,$$
since $z_1>z_2$. 

\hfill$\Box$\\ 

\begin{lem}\label{increase}
Under   \eqref{basiccoupling}, if $\rb(j)-\rb(i)\leq \tilde\rb(j)-\tilde\rb(i)$, for all $i\leq j$, then  
$$0\leq \left(\bar L^{\tilde\rb}[i]_n-\bar L^{\tilde\rb}[0]_n\right)-\left(\bar L^{\rb}[i]_n-\bar L^{\rb}[0]_n\right)\leq \left(\bar L^{\tilde\rb}[j]_n-\bar L^{\tilde\rb}[0]_n\right)-\left(\bar L^{\rb}[j]_n-\bar L^{\rb}[0]_n\right)\,,$$
for all $0\leq i\leq j$, and 
$$0\leq\left(\bar L^{\rb}[i]_n-\bar L^{\rb}[0]_n\right)-\left(\bar L^{\tilde\rb}[i]_n-\bar L^{\tilde\rb}[0]_n\right)\leq \left(\bar L^{\rb}[j]_n-\bar L^{\rb}[0]_n\right)-\left(\bar L^{\tilde\rb}[j]_n-\bar L^{\tilde\rb}[0]_n\right)\,,$$
for all $j\leq i\leq 0$. 
\end{lem}

\noindent{\bf Proof\,\,} It follows from Lemma \ref{attract} since the first inequality is equivalent to 
$$\bar L^{\rb}[j]_n-\bar L^{\rb}[i]_n\leq \bar L^{\tilde\rb}[j]_n-\bar L^{\tilde\rb}[i]_n\,,$$
for $i\leq j$, while the second one is equivalent to 
$$\bar L^{\rb}[i]_n-\bar L^{\rb}[j]_n\leq \bar L^{\tilde\rb}[i]_n-\bar L^{\tilde\rb}[j]_n\,,$$
for $j\leq i$. Notice that $i\geq 0$ implies that
$$0\leq \left(\bar L^{\tilde\rb}[i]_n-\bar L^{\tilde\rb}[0]_n\right)-\left(\bar L^{\rb}[i]_n-\bar L^{\rb}[0]_n\right)\,,$$
and $i\leq 0$ implies that 
$$0\leq\left(\bar L^{\rb}[i]_n-\bar L^{\tilde\rb}[0]_n\right)-\left(\bar L^{\tilde\rb}[i]_n-\bar L^{\tilde\rb}[0]_n\right)\,.$$

\hfill$\Box$\\ 

\subsection{KPZ Localisation of Exit Points}\label{ControlExit}

For $x\in \R$ and $t\geq 0$, we denote 
$$[x]_{t}\equiv (\lfloor t\rfloor+\lfloor x\rfloor , \lfloor t\rfloor-\lfloor x\rfloor)\,.$$  
Recall the definition of $Z^{\rh}$ given in Section \ref{MaxPath}, denote 
$$\yy^{\pm}_{t,n}=[\pm(2tn)^{2/3}]_{tn}=\Big(\lfloor tn\rfloor+\pm \lfloor(2tn)^{2/3}\rfloor,\lfloor tn\rfloor-\pm\lfloor(2tn)^{2/3}\rfloor\Big)\,,$$
and let 
$$\phi_t(r):=\limsup_{n} \P\Big(|Z^{\rh^{(n)}}(\yy^+_{t,n})|>r(tn)^{2/3}\Big)+\limsup_{n} \P\Big(|Z^{\rh^{(n)}}(\yy^-_{t,n})|>r(tn)^{2/3}\Big)\,.$$
The proof of the following proposition is postpone to Section \ref{ProofProp}.

\begin{prop}\label{InterControl2}
Under assumption \eqref{Assump}, there exist $t_0>0$ and $\phi(r)$ such that $\phi_t(r)\leq  \phi(r)$ for all $t\geq t_0$ and $\lim_{r\to\infty}\phi(r)=0$.
\end{prop}

In the stationary LPP model with boundary, the macroscopic location of the exit point $Z^{\rho}(n,n)$ is $(1-d_\rho)n$, where 
$$d_\rho = \left(\frac{1-\rho}{\rho}\right)^2\,,$$
and the order of fluctuations is $n^{2/3}$ \cite{BaCaSe}. Thus, if one sets $\rho = 1/2\pm c rn^{-1/3}$ then 
$Z^{\rho}(n,n)$ fluctuates around $(1-d_\rho)n\sim \pm r n^{2/3}$, which allows us to tune $\rho$ in such way that, with high probability, $Z^{\rho^-}(n,n)\leq -\epsilon rn^{2/3}\leq \epsilon r n^{2/3} \leq Z^{\rho^+}(n,n)$ for some small $\epsilon>0$.

\begin{lem}\label{stcontrol} 
Set $\rho_n^{\pm}:=\frac{1}{2}\pm c\frac{r}{n^{1/3}}$, where $c>0$ is fixed. There exist constants $\epsilon_1,C_1>0$ such that, for all $r>1$, 
$$\limsup_{n\to\infty}\P\left(\bar Z^{\rho_n^+}[-(2n)^{2/3}]_n< \epsilon_1 rn^{2/3} \right)\leq \frac{C_1}{r^3}\,,$$
and 
$$\limsup_{n\to\infty}\P\left(\bar Z^{\rho_n^-}[(2n)^{2/3}]_n>- \epsilon_1 rn^{2/3} \right)\leq \frac{C_1}{r^3}\,.$$ 
\end{lem}

\noindent{\bf Proof\,\,} The proof of this lemma follows the same argument in the proof of Lemma 2.3 \cite{Pi}, that is based on the fluctuations results for exit-points provided by \cite{BaCaSe}. One can improve the upper bound to $Ce^{-c r^2}$ by using Lemma 2.5 and Corollary 2.6 in \cite{FeOc}.    

\hfill$\Box$\\ 
	
By \eqref{Recover} we have that, 
\begin{eqnarray*}
\phi_t(r)&:=&\limsup_{n} \P\Big(|Z^{\rh^{(n)}}(\yy^+_{t,n})|>r(tn)^{2/3}\Big)+\limsup_{n} \P\Big(|Z^{\rh^{(n)}}(\yy^-_{t,n})|>r(tn)^{2/3}\Big)\\
&=&\limsup_{n} \P\Big(|\bar Z^{\rh^{(n)}}(\yy^+_{t,n})|>r(tn)^{2/3}\Big)+\limsup_{n} \P\Big(|\bar Z^{\rh^{(n)}}(\yy^-_{t,n})|>r(tn)^{2/3}\Big)\,.
\end{eqnarray*}

\begin{lem}\label{localisation}
For $r>0$ let 
$$\rho^{\pm}_{t,n}\equiv\rho_{t,n}^{\pm}(r):=\frac{1}{2}\pm c\frac{r}{\lfloor tn\rfloor^{1/3}}\,,$$ 
and 
\begin{equation}\label{CompEvent}
E_{t,n}(r)\equiv E_{t,n}^{\rh^{(n)}}(r):=\left\{\bar Z^{\rho_{t,n}^+}(\yy^-_{t,n})\geq \bar Z^{\rh^{(n)}}(\yy^+_{t,n})\,\,\mbox{ and }\,\, \bar Z^{\rho_{t,n}^-}(\yy^+_{t,n}) \leq \bar Z^{\rh^{(n)}}(\yy^-_{t,n})\right\}\,.
\end{equation}
Then,  
$$\limsup_{n\to\infty}\P\left(E_{t,n}(r)^c\right)\leq 2C_1 r^{-3}+\phi_t\left(2^{-1}\epsilon_1 r\right)\,.$$
where $\epsilon_1$ and $C_1$ are given by Lemma \ref{stcontrol}. 
\end{lem}

\noindent{\bf Proof\,\,}  The probability of $E_{t,n}(r)^c$ is bounded from above by 
\begin{equation}\label{loc1}
\P\left(\bar Z^{\rho_{t,n}^+}(\yy^-_{t,n})< \bar Z^{\rh^{(n)}}(\yy^+_{t,n}) \right)+\P\left(\bar Z^{\rho_{t,n}^-}(\yy^+_{t,n}) > \bar Z^{\rh^{(n)}}(\yy^-_{t,n})\right)\,.
\end{equation}
To deal with the first term in \eqref{loc1}, we note that  
\begin{eqnarray*}
\P\left(\bar Z^{\rho_{t,n}^+}(\yy^-_{t,n})< \bar Z^{\rh^{(n)}}(\yy^+_{t,n}) \right)&\leq &\P\left( |\bar Z^{\rh^{(n)}}(\yy^+_{t,n})|>2^{-1}\epsilon_1 r (tn)^{2/3}\right)\\
&+&\P\left(\bar Z^{\rho_{t,n}^+}(\yy^-_{t,n})< 2^{-1}\epsilon_1 r (tn)^{2/3}\right)\,.
\end{eqnarray*}
Taking $m=\lfloor tn\rfloor$ and using Lemma \ref{stcontrol} ($2^{-1}(nt)^{2/3}\leq \lfloor nt\rfloor^{2/3}$ for large $n$)
$$\limsup_{n\to\infty}\P\left(\bar Z^{\rho_{t,n}^+} (\yy^-_{t,n})< 2^{-1}\epsilon_1 r (tn)^{2/3} \right)\leq \limsup_{m\to\infty}\P\left(\bar Z^{\rho_{m}^+} [-(2m)^{2/3}]_{m}<\epsilon_1 r m^{2/3} \right)\,.$$
The analog upper bound for the second term in \eqref{loc1} is obtained by using that 
\begin{eqnarray*}
\P\left(\bar Z^{\rho_{t,n}^-}(\yy^+_{t,n}) > \bar Z^{\rh^{(n)}}(\yy^-_{t,n})\right)&\leq &\P\left( |\bar Z^{\rh^{(n)}}(\yy^-_{t,n})|>2^{-1}\epsilon_1 r (tn)^{2/3}\right) \\
&+& \P\left(\bar Z^{\rho_{t,n}^-}(\yy^+_{t,n})>-2^{-1}\epsilon_1 r (tn)^{2/3}\right)\,.
\end{eqnarray*}

\hfill$\Box$\\ 

\section{Proof of the Theorem \ref{Coupling}}

Now we are ready to construct the coupling $\left(\fh_t(\cdot;\fh),\fh_t(\cdot;\fb)\right)$ and get the upper bound \eqref{Coupling1}. Given $\fh\in UC$, we pick a sequence of sequence of particle configurations $\eta^{(n)}$ such that \eqref{kpz1c} holds, and set $\rh^{(n)}=\rh(\eta^{(n)})$. By \eqref{kpz3},
$$H_{t,n}\left(x\right)\equiv H_{t,n}\left(x;\rh^{(n)}\right) = \frac{L^{\rh^{(n)}}[2^{2/3} xn^{2/3}]_{\lfloor tn\rfloor}-4nt}{2^{4/3} n^{1/3}}\stackrel{dist.}{\to}\fh_t(x;\fh)\,$$
(as a continuous process in $x$). For the two-sided Brownian motion $\fb$ we pick a  particle configuration $\eta^{1/2}$ given by i.i.d. Bernoulli of parameter $1/2$, and set $\rh^{1/2}=\rh(\eta^{1/2})$. Thus, by \eqref{kpz3},
$$H^{1/2}_{t,n}\left(x\right)\equiv H_{t,n}\left(x;\rh^{1/2}\right) = \frac{L^{\rh^{1/2}}[2^{2/3} xn^{2/3}]_{\lfloor tn\rfloor}-4nt}{2^{4/3} n^{1/3}}\stackrel{dist.}{\to}\fh_t(x;\fb)\,$$
(as a continuous process in $x$). Therefore the pair $\left(H_{t,n},H^{1/2}_{t,n}\right)$ is relatively compact in $\cal C\times \cal C$, which implies the existence of a weak limit in $\cal C\times \cal C$, whose marginals are obviously given by $\fh_t(\cdot;\fh)$ and $\fh_t(\cdot;\fb)$, as described by \eqref{Int}. We will show next that such a coupling $\left(\fh_t(\cdot;\fh),\fh_t(\cdot;\fb)\right)$ satisfies \eqref{Coupling1}.
\newline 

Set $\delta_t:=at^{-2/3}$, $r=r_t:=\delta_t^{-1/4}$ and $\rho_{n,t}^{\pm}:=\rho_{n,t}^{\pm}(r_t)$ (the value of $c>0$ will be given later). Given a profile $\omega_\xx^{1/2}$ as in Example 1, define the boundary conditions $\rb^{\pm}$ by setting 
$$\omega^{\pm}_\xx=\left\{\begin{array}{lll}
\frac{1}{2\rho_{t,n}^\pm}\omega^{1/2}_{(0,|z|)} & \mbox{ for } z<0\\
\frac{1}{2(1-\rho_{t,n}^\pm)}\omega^{1/2}_{(z,0)}&\mbox{ for } z>0\,.\end{array}\right.$$
Thus, $\omega^{\pm}\stackrel{dist.}{=}\omega^{\rho^\pm_{n,t}}$. Since $\rho^-_{n,t}\leq 1/2\leq \rho^{+}_{t,n}$, if $i<0$ then,
$$\omega^{-}_{(0,-i)}-\omega^{1/2}_{(0,-i)}=\left(\frac{1}{2\rho^-_{n,t}}-1\right)\omega^{1/2}_{(0,-i)}\,\geq\, 0\,,$$
and
$$\omega^{+}_{(0,-i)}-\omega^{1/2}_{(0,-i)}=\left(\frac{1}{2\rho^+_{n,t}}-1\right)\omega^{1/2}_{(0,-i)}\,\leq\, 0\,.$$ 
If $i>0$,
$$\omega^{-}_{(i,0)}-\omega^{1/2}_{(i,0)}=\left(\frac{1}{2(1-\rho^-_{n,t})}-1\right)\omega^{1/2}_{(0,-i)}\,\leq\,0\,,$$
and
$$\omega^{+}_{(i,0)}-\omega^{1/2}_{(i,0)}=\left(1-\frac{1}{2(1-\rho^+_{n,t})}\right)\omega^{(1/2)}_{(i,0)},\geq\, 0\,.$$ 
Hence, if $i<j$ then 
\begin{equation}\label{coupling}
\rb^{-}(j)-\rb^{-}(i)\leq \rb^{+}(j)-\rb^{+}(i)\,.
\end{equation}

Again for simple notation, we use the superscript $\pm$ for quantities that are related to $\rho_{n,t}^{\pm}$, and the last-passages times $\bar L^{\pm}\equiv \bar L^{\rho_{n,t}^{\pm}}$, $\bar L^{\rh^{(n)}}=L^{\rh^{(n)}}$ and $\bar L^{\rh^{1/2}}=L^{\rh^{1/2}}$ (recall {\bf Example 1} and {\bf Example 2}) using the basic coupling \eqref{basiccoupling}. By \eqref{coupling} and Lemma \ref{increase}, for $x\in [0,a]$, 
\begin{equation}\label{AttractCoup1}
0\leq \Delta\bar H_{t,n}^{+}(x)-\Delta \bar H_{t,n}^{-}(x)\leq \Delta\bar H_{t,n}^{+}(a)-\Delta \bar H_{t,n}^{-}(a)\,,
\end{equation}
while for $x\in [-a,0]$,  
\begin{equation}\label{AttractCoup2}
0\leq \Delta\bar H_{t,n}^{-}(x)-\Delta \bar H_{t,n}^{+}(x)\leq \Delta\bar H_{t,n}^{-}(-a)-\Delta \bar H_{t,n}^{+}(-a)\,.
\end{equation}

Recall Lemma \ref{localisation} and let $E_{t,n}\equiv E_{t,n}^{\rh^{(n)}}(r_t)$ and $E_{t,n}^{1/2}\equiv E_{t,n}^{\rh^{1/2}}(r_t)$. Since $\delta_t:=at^{-2/3}$, then $an^{2/3}=\delta_t (tn)^{2/3}$ and $\delta_t\leq 1$ for $t\geq a^{3/2}$. Thus, 
$$2^{2/3}an^{2/3}\leq (2tn)^{2/3}\,\,\mbox{ and }\,\,-2^{2/3}an^{2/3}\geq -(2tn)^{2/3}\,.$$
By \eqref{IncExit}, this implies that, on the event $E_{t,n}$,
$$\bar Z^{+}[-2^{2/3}an^{2/3}]_{tn}\geq \bar Z^{\rh^{(n)}}[2^{2/3}an^{2/3}]_{tn}\,\,\mbox{ and }\,\, \bar Z^{-}[2^{2/3}an^{2/3}]_{tn} \leq \bar Z^{\rh^{(n)}}[-2^{2/3}an^{2/3}]_{tn}\,,$$
and on the event $E_{t,n}^{1/2}$, 
$$\bar Z^{+}[-2^{2/3}an^{2/3}]_{tn}\geq \bar Z^{\rh^{1/2}}[2^{2/3}an^{2/3}]_{tn}\,\,\mbox{ and }\,\, \bar Z^{-}[2^{2/3}an^{2/3}]_{tn} \leq \bar Z^{\rh^{1/2}}[-2^{2/3}an^{2/3}]_{tn}\,.$$

Now pick $\omega\in E_{t,n}$. By  \eqref{IncExit}, for every $x\in[-a,a]$,
$$\bar Z^{\rh^{(n)}}[2^{2/3}xn^{2/3}]_{tn}\leq \bar Z^{\rh^{(n)}}[2^{2/3}an^{2/3}]_{tn}\,\,,\,\, \bar Z^{\rh^{(n)}}[-2^{2/3}an^{2/3}]_{tn}\leq \bar Z^{\rh^{(n)}}[2^{2/3}xn^{2/3}]_{tn}\,,$$ 
and 
$$\bar Z^{+}[-2^{2/3}an^{2/3}]_{tn}\leq \bar Z^{+}[2^{2/3}xn^{2/3}]_{tn}\,\,,\,\, \bar Z^{-}[2^{2/3}xn^{2/3}]_{tn}\leq \bar Z^{-}[2^{2/3}an^{2/3}]_{tn}\,.$$
Hence
$$\bar Z^{\rh^{(n)}}[2^{2/3}xn^{2/3}]_{tn}\leq \bar Z^{+}[0]_{tn}\,\,\mbox{ and }\,\,\bar Z^{-}[2^{2/3}xn^{2/3}]_{tn}\leq \bar Z^{\rh^{(n)}}[0]_{tn}\,,\,\mbox{ for $x\in[0,a]$}\,,$$
while 
$$\bar Z^{\rh^{(n)}}[0]_{tn}\leq \bar Z^{+}[2^{2/3}xn^{2/3}]_{tn}\,\,\mbox{ and }\,\,\bar Z^{-}[0]_{tn}\leq \bar Z^{\rh^{(n)}}[2^{2/3}xn^{2/3}]_{tn}\,,\,\mbox{ for $x\in[-a,0]$}\,.$$
By Lemma \ref{comparison}, for $x\in[0,a]$, (recall \eqref{Recover})
$$\bar L^{-}[2^{2/3}xn^{2/3}]_{tn}-\bar L^{-}[0]_{tn}\leq  L^{\rh^{(n)}}[2^{2/3}xn^{2/3}]_{tn}- L^{\rh^{(n)}}[0]_{tn} \leq \bar L^{+}[2^{2/3}xn^{2/3}]_{tn}-\bar L^{+}[0]_{tn}\,,$$
while for $x\in[-a,0]$,
$$\bar L^{-}[0]_{tn}-\bar L^{-}[2^{2/3}xn^{2/3}]_{tn}\leq  L^{\rh^{(n)}}[0]_{tn}- L^{\rh^{(n)}}[2^{2/3}xn^{2/3}]_{tn} \leq \bar L^{+}[0]_{tn}-\bar L^{+}[2^{2/3}xn^{2/3}]_{tn}\,.$$
Therefore,
$$\Delta\bar H_{t,n}^{-}(x)\leq \Delta H_{t,n}(x)\leq \Delta\bar H_{t,n}^{+}(x)\,,\mbox{ for $x\in[0,a]$}\,,$$
and
$$\Delta\bar H_{t,n}^{+}(x)\leq \Delta H_{t,n}(x)\leq \Delta\bar H_{t,n}^{-}(x)\,,\mbox{ for $x\in[-a,0]$}\,.$$
Repeating the same argument, if we pick $\omega\in E_{t,n}^{1/2}$ we get 
$$\Delta\bar H_{t,n}^{-}(x)\leq \Delta H_{t,n}^{1/2}(x)\leq \Delta\bar H_{t,n}^{+}(x)\,,\mbox{ for $x\in[0,a]$}\,,$$
and
$$\Delta\bar H_{t,n}^{+}(x)\leq \Delta H_{t,n}^{1/2}(x)\leq \Delta\bar H_{t,n}^{-}(x)\,,\mbox{ for $x\in[-a,0]$}\,.$$
Together with \eqref{AttractCoup1} and \eqref{AttractCoup2}, this implies that, for $\omega\in E_{t,n}\cap  E_{t,n}^{1/2}$, 
\begin{eqnarray*}
|\Delta H_{t,n}(x)-\Delta H_{t,n}^{1/2}(x)|&\leq &\left(\Delta \bar H_{t,n}^{+}(x)-\Delta\bar H_{t,n}^{-}(x)\right)\mathds{1}\{x\in[0,a]\}\\
&+& \left(\Delta\bar H_{t,n}^{-}(x)-\Delta\bar H_{t,n}^{+}(x)\right)\mathds{1}\{x\in[-a,0]\}\\
&\leq &\Delta\bar H_{t,n}^{+}(a)-\Delta\bar H_{t,n}^{-}(a)\\
&+ &\Delta\bar H_{t,n}^{-}(-a)-\Delta\bar H_{t,n}^{+}(-a)\,,
\end{eqnarray*}
and hence,  
$$\sup_{x\in[-a,a]}|\Delta H_{t,n}(x)-\Delta H_{t,n}^{1/2}(x)|\leq I_{t,n}(a)\,,$$
where 
$$I_{t,n}(a):=\Delta\bar H_{t,n}^{+}(a)-\Delta\bar H_{t,n}^{-}(a)+\Delta\bar H_{t,n}^{-}(-a)-\Delta\bar H_{t,n}^{+}(-a)\,.$$
Notice that 
$$I_{t,n}(a)=\frac{\left(\bar L^{+}[2^{2/3}an^{2/3}]_{tn}-\bar L^{+}[-2^{2/3}an^{2/3}]_{tn}\right)-\left(\bar L^{-}[2^{2/3}an^{2/3}]_{tn}-\bar L^{-}[-2^{2/3}an^{2/3}]_{tn}\right)}{2^{4/3}n^{1/3}}\,,$$
and thus, by Lemma \ref{attract}, $I_{t,n}(a)\geq 0$. Therefore,
\begin{eqnarray}
\nonumber\P\left(\sup_{x\in[-a,a]}|\Delta H_{t,n}(x)-\Delta H_{t,n}^{1/2}(x)|>\eta\right)&\leq & \P\left(E_{t,n}^c\right)+\P\left((E_{t,n}^{1/2})^c\right)+\P\left(I_{t,n}(a)>\eta\right)\\
\label{markov} &\leq &\P\left(E_{t,n}^c\right)+\P\left((E_{t,n}^{1/2})^c\right)+\frac{\E\left(I_{t,n}(a)\right)}{\eta}\,.
\end{eqnarray}
 
To control $\E\left(I_{t,n}(a)\right)$, write $2^{4/3}n^{1/3}I_{t,n}(a)$ as a sum of the increments \eqref{statinc}:
$$2^{4/3}n^{1/3}I_{t,n}(a)=\sum_{k=-\lfloor 2^{2/3}an^{2/3}\rfloor+1}^{\lfloor 2^{2/3}an^{2/3}\rfloor}\left(\zeta^{+}_{k+1,n}-\zeta^{-}_{k+1,n}\right)\,.$$
By \eqref{stat}, 
\begin{eqnarray*}
\E\left(\zeta^{+}_{k,n}-\zeta^{-}_{k,n}\right)&=& \left(\frac{1}{1-\rho_{t,n}^+}-\frac{1}{1-\rho_{t,n}^-}\right)+\left(\frac{1}{\rho_{t,n}^-}-\frac{1}{\rho_{t,n}^+}\right)\\
&=& \left(\frac{1}{(1-\rho_{t,n}^+)(1-\rho_{t,n}^-)}+\frac{1}{\rho_{t,n}^+\rho_{t,n}^-}\right)\left( \rho_{t,n}^+-\rho_{t,n}^-\right)\,\,=\,\,2\frac{\rho_{t,n}^+-\rho_{t,n}^-}{\rho_{t,n}^+\rho_{t,n}^-}\,.
\end{eqnarray*}
Since $\rho_{t,n}^\pm\to 1/2$, as $n\to\infty$,
$$0\leq \E\left(\zeta^{+}_{k,n}-\zeta^{-}_{k,n}\right)\leq 9 \left( \rho_{t,n}^+-\rho_{t,n}^-\right)=18c \frac{\delta^{-1/4}_t}{(tn)^{1/3}}\,,$$
for large enough $n$, and this shows that 
\begin{equation*}
2^{4/3}n^{1/3}\E\left(I_{n,t}(a)\right)\leq \left(2\times 2^{2/3}an^{2/3}\right)\times \left(18 c\frac{\delta^{-1/4}_t}{(tn)^{1/3}}\right)\,\Rightarrow \E\left(I_{n,t}(a)\right)\leq 23 c\frac{\sqrt{a}\delta_t^{1/4}}{\eta}\,,
\end{equation*}
for large enough $n$, where we use that $at^{-1/3}\delta_t^{-1/4}=\sqrt{a}\delta_t^{1/4}$. Together with \eqref{markov}, this finally yields 
$$\P\left(\sup_{x\in[a,-a]}\Delta H_{t,n}(x)-\Delta H_{t,n}^{1/2}(x)|>\eta\right)\leq \P\left(E_{t,n}(r)^c\right)+\P\left(E_{t,n}^{1/2}(r)^c\right)+\frac{\sqrt{a}\delta_t^{1/4}}{\eta}\,,$$
for large enough $n$ and $c=23^{-1}$. By Lemma \ref{localisation} (recall that $r_t:=\delta_t^{-1/4}$), 
$$\limsup_{n}\P\left(E_{t,n}(r)^c\right)+\P\left(E_{t,n}^{1/2}(r)^c\right)\leq 4C_1 \delta_t^{3/4}+\phi_t(\delta_t^{-1/4})+\phi^{1/2}_t(\delta_t^{-1/4})\,,$$
where $\phi^{1/2}_t$ is the upper bound with respect to initial profile $\rh^{1/2}$. It converges to Brownian motion, which clearly satisfies \eqref{Assump}. Since $\left(\fh_t(\cdot;\fh),\fh_t(\cdot;\fb)\right)$ is a sub-sequential weak limit point of $\left(H_{t,n},H^{1/2}_{t,n}\right)$, then 
\begin{equation*}
\P\Big(\sup_{x\in[-a,a]}|\Delta\fh_t(x;\fh)-\Delta \fh_t(x;\fb)|>\eta\sqrt{a}\Big)\leq \theta(\delta_t)+\frac{\delta_t^{1/4}}{\eta}\,,
\end{equation*}
and Proposition \ref{InterControl2} completes the proof of Theorem \ref{Coupling}, with $\theta(\delta)=4C_1 \delta^{3/4}+\phi(\delta^{-1/4})+\phi^{1/2}(\delta^{-1/4})$.

\section{Proof of Proposition \ref{InterControl2}} \label{ProofProp}

The proof Proposition \ref{InterControl2} relies on tail estimates for the location of the maxima in the definition of the last-passage percolation time. Write $\rh_\dd(k):=\rh(k)+(1,1)$ and 
$$L^{\rh}(\xx):=\max_{k\in\ZZ\,:\,\rh_\dd(k)\leq \xx} L\left(\rh_\dd(k),\xx\right)\,.$$
Let $K^{\pm}_{t,n}\in\ZZ$ be the index of the location of the maxima for $\yy^{\pm}_{t,n}$, i.e. 
$$L^{\rh}\left(\yy^{\pm}_{t,n}\right)=L\left(\rh_\dd\left(K^{\pm}_{t,n}\right)\,,\, \yy^{\pm}_{t,n}\right)\,.$$ 
By definition, $Z^{\rh^{(n)}}(\yy^\pm_{t,n})$ belongs to the intersection between the maximal path $\pi^{\rh^{(n)}}\left(\yy^\pm_{t,n}\right)$, starting at $\rh_\dd\left(K^{\pm}_{t,n}\right)$ and ending at $\yy^\pm_{t,n}$, and the non-negative coordinate axis. If we are able to control the order of magnitude of $K^{\pm}_{t,n}$ then Proposition \ref{InterControl2} will follow from well known results on the fluctuations of point to point maximal paths. 
\begin{lem}\label{exitcontrolP2}   
Under \eqref{Assump} there exist $\psi_1(r)$ and $t_0>0$ such that for all $t>t_0$  
$$\limsup_{n\to\infty}\P\left(|K^{\pm}_{n,t}|\geq r(tn)^{2/3} \right)\leq \psi_1(r)\,\,\mbox{ and }\,\,\lim_{r\to\infty}\psi_1(r)=0\,.$$
\end{lem} 

\begin{figure}{h}
\centering
\includegraphics[height=8cm]{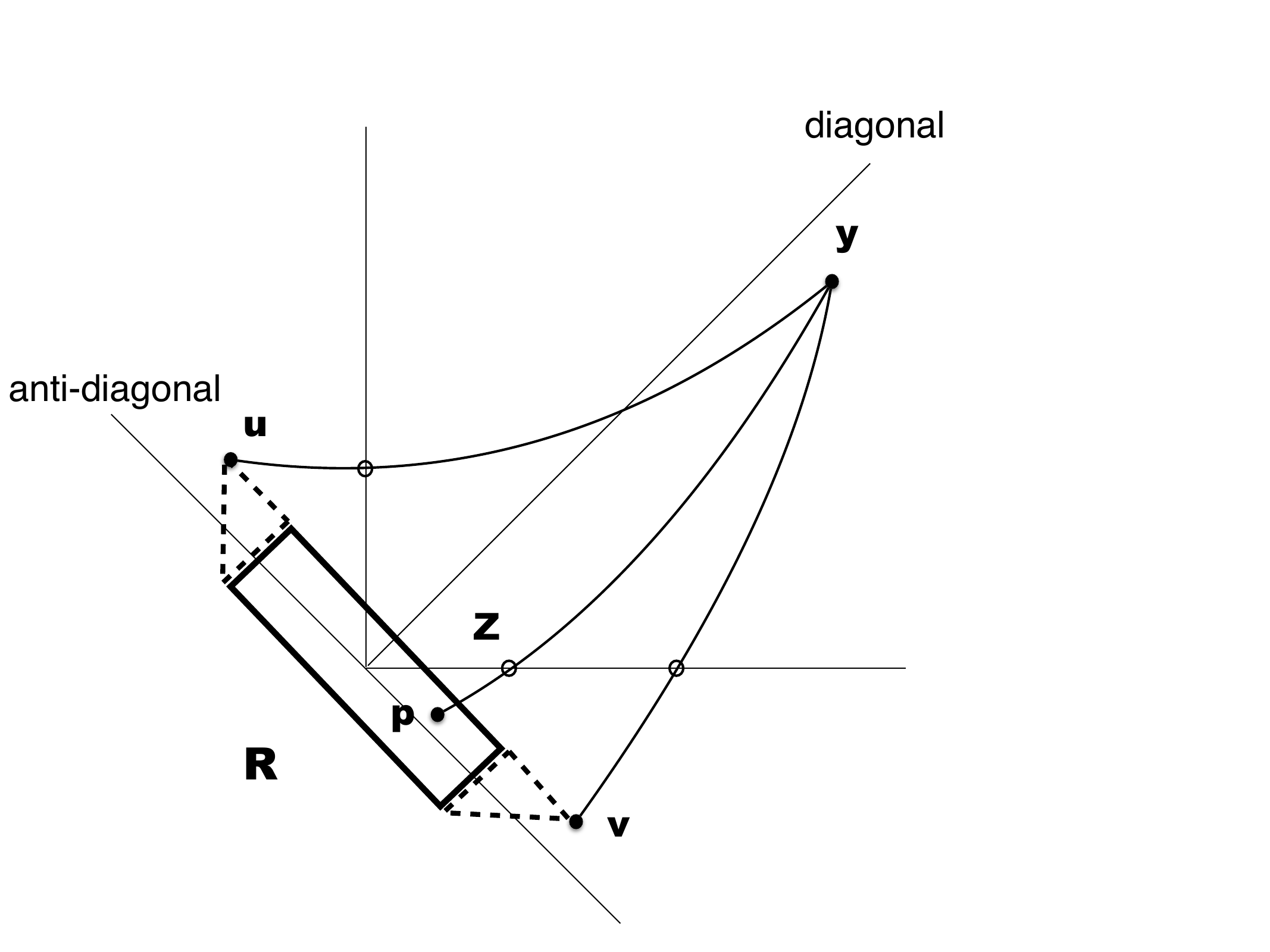}
\caption{Enclosing a maximal path starting at some point $\pp\in\mathbf{R}$ by the  point-to-point maximal paths starting at $\uu$ and $\vv$.}\label{Geo}
\end{figure}

Let $\mathbf{R}=\mathbf{R}(n,t,r)\subseteq \R^2$ be a rectangle centered at the origin and of size of $r(tn)^{2/3}$ in the anti-diagonal direction and of size $r(tn)^{1/3}$ in the diagonal direction (Figure \ref{Geo}).
By combining Lemma \ref{exitcontrolP2} together with \eqref{kpz1c}, one has that there exist $t_0>0$ and $\psi_2(r)$ such that for all $t\geq t_0$
\begin{equation}\label{EncGeo}
\limsup_{n\to\infty}\P\left(\rh_\dd\left(K^{\pm}_{t,n}\right)\in \mathbf{R}\right)\leq \psi_2(r)\,,\,\mbox{ and }\,\lim_{r\to\infty}\psi_2(r)=0\,.
\end{equation}
To conclude the proof of Proposition \ref{InterControl2}, we need to show that for all point $\pp\in\mathbf{R}$, that is a corner of $\rh$ differently from $(1,1)$, the first intersection (following the up-right orientation) between the point-to-point maximal path $\pi(\pp,\yy)$ ($\yy\equiv\yy^{\pm}_{t,n}$) and the non-negative coordinate axis is at distance of order $(tn)^{2/3}$ from the origin.  Now, let $\uu$ and $\vv$ be as it is indicated in Figure \ref{Geo}. Notice that both $\uu$ and $\vv$ are at a distance of order $(tn)^{2/3}$ from the origin, and that every maximal path starting at some point $\pp\in\mathbf{R}$ and ending at $\yy$ is below $\pi(\uu,\yy)$ and above $\pi(\vv,\yy)$. Therefore, the intersection of $\pi(\pp,\yy)$ with the non-negative coordinate axis lyes in between the respective intersections of $\pi(\uu,\yy)$ and $\pi(\vv,\yy)$. To prove that these intersections are at a distance of order $(tn)^{2/3}$ from the origin one only needs to use upper bounds for fluctuations of point-to-point maximal paths (Theorem 2.5 \cite{BaCaSe}).       
\newline 
    
\noindent{\bf Proof of Lemma \ref{exitcontrolP2}\,\,} Define 
$$L_r^{\rh}(\xx):=\max_{|k|>rn^{2/3}:\,\rh_\dd(k)\leq \xx} L\left(\rh_\dd(k),\xx\right)\,,$$ 
$$H^r_{t,n}\left(x\right) = \frac{L_r^{\rh^{(n)}}[2^{2/3} xn^{2/3}]_{\lfloor tn\rfloor}-4nt}{2^{4/3} n^{1/3}}\,,$$
and 
$$H^{\rnw}_{t,n}\left(x\right)= \frac{ L\left(\0,[2^{2/3} xn^{2/3}]_{\lfloor tn\rfloor}\right)-4nt}{2^{4/3} n^{1/3}}\,.$$
Thus
\begin{eqnarray}
\nonumber\P\left(|K^{\pm}_{t,n}| >rt^{2/3}n^{2/3}\right)&\leq &\P\left(L\left(\0,\yy^\pm_{t,n}\right)<L_{rt^{2/3}}^{\rh^{(n)}}\left( \yy^\pm_{t,n}\right)\right)\\
\label{Prop1&2}&=&\P\left(H^{\rnw}_{t,n}\left(\pm t^{2/3}\right)< H^{rt^{2/3}}_{t,n}\left(\pm t^{2/3}\right)\right)\,,
\end{eqnarray}
since $L\left(\0,\yy^\pm_{t,n}\right)\leq L^{\rh^{(n)}}\left(\yy^\pm_{t,n}\right)$ and  $L^{\rh}\left(\yy^\pm_{t,n}\right)=L_r^{\rh^{(n)}}\left(\yy^\pm_{t,n}\right)$ if  $|K^{\pm}_{t,n}| >r(tn)^{2/3}$. The LPP time $L\left(\0,[2^{2/3} xn^{2/3}]_{\lfloor tn\rfloor}\right)$ corresponds to the narrow wedge initial profile and    
$$\lim_{n\to\infty}H^{\rnw}_{t,n}\left(x\right)\stackrel{dist.}{=}t^{1/3}A(xt^{-2/3})-\frac{x^2}{t}\,,$$
while for $L_{rt^{2/3}}^{\rh^{(n)}}\left( \yy^\pm_{t,n}\right)$ we have  (by \eqref{VarFor})
$$\lim_{n\to\infty}H^{rt^{2/3}}_{t,n}\left(x\right)\stackrel{dist.}{=}\fh^{rt^{2/3}}_t(x;\fh)\stackrel{dist.}{=}\sup_{|z|>rt^{2/3}}\left\{\fh(z)+t^{1/3}A(zt^{-2/3})-\frac{(x-z)^2}{t}\right\}\,.$$
Since for any $R>0$,
$$\P\left(H^{\rnw}_{t,n}\left(\pm t^{2/3}\right)< H^{rt^{2/3}}_{t,n}\left(\pm t^{2/3}\right)\right)\leq \P\left(H^{\rnw}_{t,n}\left(\pm t^{2/3}\right)\leq -R\right)+\P\left( -R< H^{rt^{2/3}}_{t,n}\left(\pm t^{2/3}\right)\right)\,,$$
by \eqref{Prop1&2}, we get that
\begin{equation}\label{1Prop1&2}
\limsup_{n}\P\left(|K^{\pm}_{t,n}| >r(tn)^{2/3}\right)\leq  \P\left( t^{1/3}\left(A(0)-1\right)\leq -R\right)+\P\left( -R< \fh^{rt^{2/3}}_t\left(\pm t^{2/3};\fh\right)\right)\,.
\end{equation}

By using \eqref{1Prop1&2} with $R=t^{1/3} r^2/8$ we have that 
$$\limsup_{n}\P\left(|K^{\pm}_{t,n}| >r(tn)^{2/3}\right)\leq  \P\left(A(0)\leq -\frac{r^2}{8}+1\right)+\P\left( -t^{1/3}\frac{r^2}{8}< \fh^{rt^{2/3}}_t\left(\pm t^{2/3};\fh\right)\right)\,.$$
Since 
$$\P\left(A(0)\leq -\frac{r^2}{8}+1\right)\to 0\,,\mbox{ as }r\to \infty\,,$$
we only need to show that
$$\P\left( -t^{1/3}\frac{r^2}{8}< \fh^{rt^{2/3}}_t\left(\pm t^{2/3};\fh\right)\right) \to 0 \,,\mbox{ as }r\to \infty\,$$
uniformly on $t$. Indeed, for $\gamma_t=t^{1/3}$, we note that 
\begin{eqnarray*}
\fh^{rt^{2/3}}_t(\pm t^{2/3};\fh)&\stackrel{dist.}{=}&\sup_{|z|>rt^{2/3}}\left\{\fh(z)+t^{1/3}A(zt^{-2/3})-\frac{(z-\pm t^{2/3})^2}{t}\right\}\\
&=&t^{1/3}\sup_{|z|>rt^{2/3}}\left\{t^{-1/3}\fh(z)+A(zt^{-2/3})-(zt^{-2/3}-\pm 1)^2\right\}\\
&=&t^{1/3}\sup_{|u|>r}\left\{S_{\gamma_t}\fh(u)+A(u)-(u-\pm 1)^2\right\}\\
&\leq&t^{1/3}\sup_{|u|>r}\left\{S_{\gamma_t}\fh(u)+A(u)-\frac{u^2}{2}\right\}\,,\\
\end{eqnarray*}
for large enough $r>0$. But if $S_{\gamma_t}\fh(u)\leq \frac{r}{8}|u|$ then $S_{\gamma_t}\fh(u)\leq \frac{u^2}{4}$ for $|u|>r>1$, and hence
$$\P\left(-t^{1/3}\frac{r^2}{8}< \fh^{rt^{2/3}}_t(\pm t^{2/3};\fh)\right)\leq \psi(r/8)+ \P\left(-\frac{r^2}{8}<\sup_{|u|>r}\left\{A(u)-\frac{u^{2}}{4}\right\} \right)\,.$$
By (b)-Proposition 2.13 \cite{CoLiWa}, 
$$\P\left(-\frac{r^2}{8}<\sup_{|u|>r}\left\{A(u)-\frac{u^{2}}{4}\right\} \right)\to 0\,,\mbox{ as }r\to\infty\,.$$

\hfill$\Box$\\


\begin{thebibliography}{10}

\bibitem{AmCoQu}
\textsc{G. Amir, I. Corwin and J. Quastel.} Probability distribution of the free energy of the continuum directed random polymer in 1 + 1 dimensions. \newblock\emph{Commun. Pure Appl. Math.} {\bf 64} (2011), 466--537.

\bibitem{BaCaSe}
\textsc{M. Bal\'azs, E.~A. Cator, and T. Sepp\"al\"ainen}. Cube root fluctuations for the corner growth model associated to the exclusion process.
\newblock \emph{Elect. J. Probab.} {\bf 11} (2006), 1094--1132.

\bibitem{Bi}
\textsc{P. Billingsley}. Convergence of probability measures.
\newblock John Wiley \& Sons, New York (1968).

\bibitem{BoCoFeVe}
\textsc{A. Borodin, I. Corwin, P. ~L. Ferrari and B. Vet\"o}. Height Fluctuations for the Stationary KPZ Equation.
\newblock \emph{Math. Phys. Anal. Geom.} {\bf 18} (2015), 20.

\bibitem{BoFePrSa}
\textsc{A. Borodin, P.~L. Ferrari, M. Pr\"ahofer and T. Sasamoto}. Fluctuation properties of the TASEP with periodic initial configuration. 
\newblock \emph{J. Statist. Phys.} {\bf 129} (2007), 1055--1080.

\bibitem{CaPi}
\textsc{E. ~A. Cator and L. ~P. ~R. Pimentel}. On the local fluctuations of last-passage percolation models.
\newblock \emph{Stoc. Proc. App.} {\bf 125} (2016), 538--551.

\bibitem{CoHa}
\textsc{I. Corwin and A. Hammond}. Brownian Gibbs property for Airy line ensembles. \newblock\emph{Invent. Math.}  {\bf 195} (2014), 441--508.

\bibitem{CoLiWa}
\textsc{I. Corwin, Z. Liu and D. Wang}. Fluctuations of TASEP and LPP with general initial data. \newblock\emph{Ann. Appl. Probab.} {\bf 26} (2016), 2030--2082.

\bibitem{CoQuRe}
\textsc{I. Corwin, J. Quastel and D. Remenik}. Renormalization fixed point of the KPZ universality class. \newblock\emph{J. Stat. Phys.} {\bf 160} (2015), 815--834.

\bibitem{DaOrVi}
\textsc{D. Dauvergne, J. Ortmann and B. Vir\'ag}. The directed landscape. \newblock Available from arXiv:1812.00309.

\bibitem{FeOc}
\textsc{P.~L. Ferrari and A. Occelli}. Universality of the GOE Tracy-Widom distribution for TASEP with arbitrary particle density. \newblock \emph{Elect. J. Probab.} {\bf 23}  (2018), 1--24.

\bibitem{Hag}
\textsc{J. H\"agg}. Local fluctuations in the Airy and discrete PNG process.
\newblock \emph{Ann. Probab.} {\bf 36} (2008), 1059--1092.

\bibitem{Ha}
\textsc{A. Hammond}. A patchwork quilt sewn from Brownian fabric: Regularity of polymer weight profiles in Brownian last passage percolation. \newblock\emph{Forum of Mathematics, Pi} {\bf 7} (2019).

\bibitem{Jo1}
\textsc{K. Johansson}. Shape fluctuations and random matrices. \newblock \emph{Comm. Math. Phys.} {\bf 209} (2000), 437--476.

\bibitem{Jo2}
\textsc{K. Johansson}. Discrete Polynuclear Growth and Determinantal processes. \newblock \emph{Comm. Math. Phys.} {\bf 242} (2003), 277--239.


\bibitem{KPZ}
\textsc{M. Kardar, G. Parisi, Y.~-C. Zhang}. Dynamic scaling of growing interfaces. \newblock\emph{Phys. Rev. Lett.} {\bf 56} (1986), 889--892.

\bibitem{MaQuRe}
\textsc{K. Matetski, J. Quastel, D. Remenik}. The KPZ fixed point.
\newblock Available from arXiv:1701.00018.

\bibitem{Pi}
\textsc{L.~P.~R. Pimentel}. Local Behavior of Airy Processes.
\newblock \newblock\emph{J. Stat. Phys.} {\bf 173} (2018), 1614--1638.

\bibitem{Sa}
\textsc{T. Sasamoto}. Spatial correlations of the 1D KPZ surface on a flat substrate. \newblock{J. Phys. A: Math. Gen.} {\bf 38} (2005), 549.


\end{thebibliography}
\end{document}